\documentclass[12pt]{amsart}
\usepackage{amscd}
\usepackage{amssymb}
\usepackage{amsmath}
\usepackage{amsthm}

\newtheorem*{theorem}{ \bf Theorem}
\newtheorem{df}{ \sc Definition}[section]

\newtheorem*{akn}{ \sc Aknowledgements}

\newtheorem{lem}[df]{ \sc Lemma}

\def\dim{{\rm dim}}
\def\deg{{\rm deg}}

\def\rank{{\rm rank}}

\def\mpr#1{\;\smash{\mathop{\hbox to 20pt{\rightarrowfill}}\limits^{#1}}\;}
\def\epi#1{\;\smash{\mathop{\hbox to 20pt{\rightarrowfill}\hskip
-15pt\rightarrow}\limits^{#1\,}}\;}

\def\Oc{{\mathcal O}}

\def\P{{\mathbb P}}

\def\Hom{{\mathcal H om}}
\def\banica{{B\u anic\u a}}

\def\qed{\hfill$\Box$\vskip10pt}

\begin{document}
\title{Globally Generated Vector Bundles on ${\P}^n$ with $c_1=3$}
\author{Cristian Anghel and Nicolae Manolache}
\address{Institute of Mathematics of the Romanian Academy
 {\em Simion Stoilow }
P.O. Box 1-764, 014700 Bucharest, Romania}
\date{}

\maketitle
\begin{abstract}
One classifies the globally generated vector bundles on $\P^n$  whith the first 
Chern class $c_1=3$ (in this paper $n\not =3$, the case $n=3$ being studied in 
\cite{Ma4}). The case $c_1=1$ is very easy, the case $c_1=2$ was done in \cite{SU}, 
the case $c_1=3$, rank$=2$ was settled in \cite{Huh} and the case $c_1 \le 5$, rank=2
in \cite{CE}.
Our work is based on Serre's theorem relating vector bundles of $\rank =2$ with 
codimension $2$ lci subschemes and its generalization 
for higher ranks, considered firstly by Vogelaar in \cite{Vog}.
\end{abstract}

\section{Introduction}
\label{intro}
We use freely the standard notation, namely that  introduced in \cite{FAC},
used (and extended  in \cite{EGA}) (cf. also \cite{HAG}).
Our base field is $\mathbb C$. 

A globally generated nontrivial vector bundle $E$ on a a projective space has the 
first Chern class $c_1(E) > 0$. When $c_1(E)=1$, modulo a trivial summand, 
$E\cong \Oc (1)$ or $E \cong T(-1)$, $T$ being the tangent bundle. The case 
$c_1(E)=2$ is studied in \cite{SU} and the case $c_1=3$ on $\P^3$ was studied in
\cite{Ma4}.  Globally generated rank $2$ vector bundles with 
$c_1=3$ are listed in \cite{Huh}.  Rank $2$ globally generated vector bundles on 
$\P^n$, $n\ge 3$, with $c_1 \le 5$  are listed in \cite{CE}.

The aim of this paper is to prove:
\begin{theorem} (i) On $\P^2$ the globally generated vector bundles of rank 2
with the first Chern class $c_1=3$, which are not direct sums of line bundles
are as follows:
\begin{align*}
\tag{1} 0 \to \Oc(-1) \to 2\Oc\oplus \Oc(2) \to & E\to 0 \\
\tag{2} E \cong \Omega ^1(3) \cong T_{\P^2}\\
\tag{3} 0 \to \Oc(-1) \to \Oc \oplus 2\Oc(1) \to & E \to 0 \\
\tag{4} 0 \to 2\Oc(-1) \to 3\Oc\oplus \Oc(1) \to & E \to 0\\
\tag{5} 0 \to \Oc(-2) \to 2\Oc\oplus \Oc(1) \to & E \to 0\\
\tag{6} 0 \to \Oc(-3) \to 3\Oc \to & E \to 0\\
\tag{7} 0 \to \Oc(-2)\oplus \Oc(-1) \to 4\Oc \to & E \to 0\\
\tag{8} 0 \to 3\Oc(-1) \to 5\Oc \to & E \to 0\\
\end{align*}

(ii) On $\P^4$ the globally generated vector bundles with $c_1=3$ of rank at most 
$4$ are:

-- none of rank $2$

-- of rank $3$: the dual of the Trautmann-Vetter (Tango)  bundle:
$$
0 \to T(-2) \to 7\Oc \to E \to 0
$$
-- of rank $4$: given by exact sequences:
\begin{align*}
\tag{1} 0 \to \Oc(-1) \to 4\Oc \oplus \Oc(2) \to  & E \to 0 \\
\tag{2}  0 \to \Oc(-1) \to 3\Oc \oplus 2\Oc(1) \to & E \to 0 \\
\tag{3}  0 \to 2\Oc(-1) \to 5\Oc \oplus \Oc(1) \to & E \to 0 \\
\tag{4}  0 \to \Oc(-2) \to 4\Oc \oplus \Oc(1) \to & E \to 0 \\
\tag{5}  0 \to 3\Oc(-1) \to 7\Oc  \to & E \to 0 \\
\tag{6}  0 \to \Oc(-2) \oplus \Oc(-1) \to 6\Oc  \to & E \to 0 \\
\tag{7}  0 \to \Oc(-3) \to 5\Oc \to & E \to 0 \\
\tag{8}  E \ \cong & \Omega ^1(2)
\end{align*}

(iii) On $\P^5$ the globally generated vector bundles with $c_1=3$ of rank at most 
$5$ are:

-- none of ranks $2$, $3$, $4$

-- of rank $5$ given by exact sequences:
\begin{align*}
\tag{1} 0 \to \Oc(-1) \to 5\Oc \oplus \Oc(2) \to & E\to 0 \\
\tag{2} 0 \to \Oc(-1) \to 4\Oc\oplus 2\Oc(1) \to & E \to 0\\
\tag{3} 0 \to 2\Oc(-1) \to 6\Oc \oplus \Oc(1) \to & E \to 0 \\
\tag{4} 0 \to \Oc(-2) \to 5\Oc\oplus \Oc(1) \to & E \to 0\\
\tag{5} 0 \to 3\Oc(-1)\to  8\Oc \to & E \to 0\\
\tag{6} 0 \to \Oc(-2) \oplus \Oc(-1) \to 7\Oc \to & E \to 0\\
\tag{7} 0 \to \Oc(-3) \to 6\Oc \to & E \to 0\\
\end{align*}

(iv) On $\P^n$, $n=2,4,5$ the globally generated vector bundles of rank $ >n$
with $c_1=3$ are extensions of those of rank $n$  by trivial bundles.

(v) On $\P^n$, $n \ge 5$ the globally generated vector bundles with $c_1=3$ 
have rank  $r \ge n$ and are given by exact sequences:
\begin{align*}
\tag{1} 0 \to \Oc(-1) \to r\Oc \oplus \Oc(2) \to & E\to 0 \\
\tag{2} 0 \to \Oc(-1) \to (r-1)\Oc\oplus 2\Oc(1) \to & E \to 0\\
\tag{3} 0 \to 2\Oc(-1) \to (r+1)\Oc \oplus \Oc(1) \to & E \to 0 \\
\tag{4} 0 \to \Oc(-2) \to r\Oc\oplus \Oc(1) \to & E \to 0\\
\tag{5} 0 \to 3\Oc(-1)\to  (r+3)\Oc \to & E \to 0\\
\tag{6} 0 \to \Oc(-2) \oplus \Oc(-1) \to (r+2)\Oc \to & E \to 0\\
\tag{7} 0 \to \Oc(-3) \to (r+1)\Oc \to & E \to 0\\
\end{align*}
\end{theorem}

\begin{akn}
The starting point of this paper was a discussion with Iustin Coand\u a about
\cite{SU}. We want to thank him for this, for his interest, for pointing us
some useful references and also for useful comments on a previous version of this paper. 
The second author expresses his thanks to the Institute of Mathematics, Oldenburg 
University, especially to Udo Vetter, for warm hospitality. In writing this paper 
the authors had no other support besides  the membership to the Institute of 
Mathematics ``Simion Stolow'' of the Romanian Academy. 
\end{akn}

\section{Preliminaries}
\label{pre}
We recall some general results, mostly of which are the starting remarks in 
\cite{SU}:\smallskip

\begin{itemize}
\item[$\mathbf A.$] (J.-P. Serre, cf. \cite{At}) If a vector bundle $E$ or rank 
$r > n$ on a smooth variety $X$ of dimension $n$ is globally generated, then $E$ 
contains a  trivial subbundle of rank $=r-n$ .

\item[$\mathbf B.$] If a  globally generated vector bundle $E$ or rank $r\le n$ 
on ${\P}^n$ has  $c_r=0$,  then $E$ contains a  trivial subbundle of rank $=1$ 
(cf. \cite{OSS}, chap. I, Lemma 4.3.2).

\item[$\bf C.$] If $E$ is a globally generated vector bundle on ${\P}^n$, then 
$E$ contains a trivial direct summand of rank $=\dim H^0(E^\vee )$ (cf. \cite{Ott}, 
Lemma 3.9).

\item[$\mathbf D.$] If $E$ is a vector bundle on $\P^n$, which is globally generated
and for which $H^0(E(-c_1))\not =0$, where $c_1$ is the first Chern class of $E$,
 then, modulo a trivial summand, $E\cong \Oc (c_1)$ (cf. \cite{S}, Prop.1).

\item[$\bf E.$] (Vogelaar, Arrondo's setting, cf. \cite{Arr}, Theorem 1.1) 

(i) If $E$ is a vector bundle of rank $r$ on $\P^n$ and $r-1$ sections 
$s_1,\ldots ,s_{r-1}$ of $E$  has as dependency locus a subscheme $Y$ of $\P^n$ 
which is lci, then one has an exact sequence:
$$
0 \to (r-1)\Oc  \xrightarrow{(s_1,\ldots ,s_{r-1})} E \to I_Y \otimes L \to 0 \ \ ,
$$
where $L=\Lambda ^rE$.
In this case, if $N$ is the normal sheaf of $Y$, then  
$\Lambda ^2 N \otimes L^\vee |_Y$ is generated by $r-1$ global sections.

(ii) Conversely, if $Y$ is a lci subscheme of codimension $2$ of $\P^n$ and $L$ is 
a line bundle on $\P^n$ with the property $H^2(L^\vee)=0$, if 
$\Lambda ^2 N \otimes L^\vee |_Y$ is generated by $r-1$ global sections, then there 
exists a vector bundle $E$ of rank $r$ such that $\Lambda ^rE=L$ and $E$ has $r-1$
global sections whose dependence locus is $Y$.

\item[$\mathbf F.$] If $\varphi : F \to G$ is a general morphism between vector 
bundles over a smooth variety and $\Hom (F,G)$ is globally generated, then the 
dependency locus $Y$ of $\varphi $ is either empty or has codimension $\rank (G) - 
\rank (F) +1$ and is nonsingular outside a subset of codimension $ \ge \rank (G) - 
\rank (F) +3$ (cf. \cite{B} 4.1 or \cite{Ch} for a more general result and 
\cite{PS}, \cite{Sa} for the case of arithmetically Cohen-Macaulay curves in $\P^3$).

\end{itemize}
\qed

According to {\bf A}, any vector bundle on $\P^n$ of rank $ >n$ is an extension of 
a rank $n$ vector bundle by a trivial one, so on $\P^n$, modulo extensions by
trivial bundles, it is enough to find the globally generated vector bundles of rank
at most $n$ with $c_1=3$.

\section{Vector Bundles on $\P^2$}
\label{Srank2}

In this case we have to treat only the case of rank =2. This is treated in 
\cite{Huh} and we give here other arguments.
Take  $E$  a globally generated vector bundle of rank 2  on $\P^2$, non-splittable 
in line bundles, with $c_1(E)=3$; we may assume $H^0(E(-3))=0$, according to 
$\mathbf D$. If no section of it has an empty zero's set (in which case E splits), 
a general section of it has a smooth zero set $Z$ which according to \cite{HM}, Prop.
5.1, or \cite{Ha}, Prop. 1.4. is a union of reduced closed points. One has an exact 
sequence:
\begin{equation}
\tag{*} 0 \to \Oc \to E \to I_Z(3) \to 0
\end{equation}

I. {\bf Assume} $H^0(E(-3))=0$, $H^0(E(-2))\not =0$. The exact sequence (*) shows
that $Z$ lies on a line. As $I_Z(3)$ is generated by global sections, $Z$ is included 
in a complete intersection of type $(1,3)$, so $c_2(E)= \deg (Z) \le 3$. 
On the other hand, as $E$ does not split, one has an exact sequence:
\begin{equation*}
 0 \to \Oc \to E(-2) \to I_X (-1) \to 0 \ \ ,
\end{equation*}
where $X$ is lci subscheme of dimension $0$ of $\P^2$. One has $c_2(E(-2))=c_2(E)-2=
\deg (X) \ge 1$. It follows $c_2(E)=3$ so that $Z$ is a complete intersection
of type $(1,3)$ and one gets for $E$ a resolution of the shape:
$$
0 \to \Oc (-1) \to 2\Oc \oplus \Oc (2) \to E \to 0 \ \ .
$$

II. {\bf Assume} $H^0(E(-2))=0$, $H^0(E(-1))\not =0$. In this case $E$ is a stable 
vector bundle of rank $2$ on $\P^2$. Let $F=E(-2)$. Then $c_1(F)=-1$ and, because of
stability,  $c_2(F) \ge 1$. If $c_2(F)=n$, 
then $F$ is given by a monad (cf. \cite{Hu}):
$$
(n-1)\Oc (-1)  \xrightarrow{a} n\Omega (1) \xrightarrow{a^T(-1)} (n-1)\Oc \ \ .
$$ 
When $n=1$, $E\cong \Omega (3)\cong T_{\P^2}$. From now on $n\ge 2$.

As $c_1(E(-1))=1$, we have an exact sequence:
\begin{equation}
 \tag{**} 0 \to \Oc \to E(-1) \to I_X(1)\to 0 \ \ ,
\end{equation}
where $X$ is a subscheme of dimension $0$ of $\deg (X)= c_2(E(-1)) =n$.

In the exact sequence (*) the scheme $Z$ consists of $c_2(E)=n+2$ points. The condition
that $Z$ lies in a complete intersection of type $(2,3)$ gives $n+2 \le 6$, so 
$2\le n \le 4$. 

If $n=2$, then $Z$ is a complete intersection of type $(2,2)$, and one gets an 
exact sequence:
$$
0 \to \Oc (-1) \to \Oc \oplus 2\Oc (1) \to E \to 0 \ \ .
$$

If n=3 then $X$ cannot lie on a line. Indeed, $h^0(I_X(1) \not =0$ would give,
via (**), $X$ being a complete intersection of type $(1,3)$, $h^0(E(-1))=2$, 
$h^0(E)=6$. From(*) would follow that $Z$, a set of 5 different points, would lie
on two distinct conics, contradiction. It follows that $X$ has a resolution of
the shape:
$$
0 \to 2\Oc (-3) \to 3\Oc (-2) \to I_X \to 0 \ \ ,
$$ 
where from $E$ has a resolution:
$$
0 \to 2\Oc(-1) \to 3\Oc \oplus \Oc (1) \to E \to 0
$$

\textsc{Remark.}  The same result could be obtained using the resolution of $I_Z$, $Z$
 being a set of $5$ points not situated on a line:
$$
0 \to 2\Oc(-4) \to 2\Oc(-3)\oplus \Oc(-2)\to I_Z \to 0 \ \ .
$$

If $n=4$, then $Z$ is a complete intersection of type $(2,3)$ and for $E$ follows
a resolution of the shape:
$$
0 \to \Oc (-2) \to 2\Oc \oplus \Oc(1) \to E \to 0 \ \ .
$$

III. {\bf Assume} $H^0(E(-1))=0$. With the notation from above, if $c_1(F)=1$, 
$c_2(F)=n \ge 2$ then $\deg (Z)=n+2 \le 9$, $Z$ being included in a complete 
intersection of type $(3,3)$, and so $n \le 7$. As $Z$ lies on no conic, 
$4 \le n \le 7$.

For $n= 4,5,6$ i.e. $\deg (Z)=6,7,8$ and the resolutions of $I_Z$ are:

(a) $deg (Z) =6$:
$$
0 \to 3\Oc(-4) \to 4\Oc(-3) \to I_Z \to 0 
$$
and then the resolution of $E$ is:
$$
0 \to 3\Oc(-1) \to 5\Oc \to E \to 0  \ \ .
$$

(b) $\deg (Z) = 7$:
$$
0 \to \Oc(-5) \oplus \Oc(-4) \to 3\Oc(-3) \to I_Z \to 0 
$$
and then the resolution of $E$ is:
$$
0 \to \Oc(-2)\oplus \Oc (-1) \to 4\Oc \to E \to 0 \ \ ;
$$

(c) $deg (Z) =8$:
$$
0 \to 2\Oc(-5) \to \Oc(-4) \oplus 2\Oc (-3) \to I_Z \to 0 
$$
and then the resolution of $E$ is:
$$
0 \to 2\Oc(-2)\to  \Oc (-1) \oplus 3\Oc \to E \to 0 \ \ .
$$

In the last case $E$ is not globally generated.

To prove (a), use the fact that  $Z$ is linked in a complete intersection of type 
$(3,3)$ to a scheme of degree $3$. In general this can be: three non-collinear 
simple points, three collinear simple points (here are included the cases when two 
or all three coincide), a double point and a simple one but not on a line, the first 
infinitesimal neighbourhood of a point. Applying \cite{PS}, 4.1, one gets the 
resolution of $I_Z$  from the resolution of the linked scheme. In our case, when 
the linked scheme is contained in a line, $Z$ would be contained in a conic, what
we excluded. In the other cases one obtains the asserted shape. The same result can 
be obtained from \cite{GH}, where more complicated cases (6 fat points!) are studied.

In the same way one proves the cases (b), (c).

If $deg (Z) =9$, then $Z$ is a complete intersection of two cubics
and the resolution of $E$ is:
$$
0 \to \Oc (-3) \to 3\Oc \to E \to 0 \ \ ,
$$

\section{Vector Bundles on $\P ^4$ or $\P^5$}
\label{Srank3}

We do not have indecomposable rank $2$ vector bundles
with $c_1=3$ on $\P^3$ (cf. \cite{Ma4} or \cite{CE}) so that, applying Horrocks' 
criterion (cf. \cite{Ho}, or \cite{OSS}, 2.3.2) they do not exist on 
$\P^n$, $n \ge 3$, either. 
In this section $\P^n$ is $\P^4$ or $\P^5$. We use {\bf E.} and {\bf F.} If $E$ is a 
globally generated vector bundle of rank $r$ on $\P^n$, with $c_1(E)=3$ and 
if $E$ has no subbundle of $\rank$ $1$, then one has an exact sequence:
\begin{equation}
\tag{***} 0 \to (r-1)\Oc \to E \to I_Z(3) \to 0 \ \ ,
\end{equation}
where $Z$ is a smooth subvariey of $\P^n$.
Observe $\Lambda _Z:=
\Lambda ^2N\otimes L^\vee |Z \cong \omega _Z \otimes \omega^{-1} _{\P^n} 
\otimes L^ \vee |Z \cong \omega _Z(n-2)$. 
Let $E$ be  non-splittable.

I. {\bf Assume} $H^0(E(-3))=0$, $H^0(E(-2))\not =0$. If $E(-2)$ has a section without
zeros, then $E$ is an extension of vector bundles:
$$
0 \to \Oc (2) \to E \to E'\to 0
$$
where $E'$ is globally generated and has $c_1=1$. Then, modulo a trivial summand, 
$E'$ is $\Oc(1)$ or $T(-1)$. In the first case $E$ is splittable, in the second has 
$\rank > n$. Suppose now that $E$ has no rank $1$ subbundle. According to {\bf E.}
and {\bf F.}, one has an extension:
$$
0 \to (r-1)\Oc \to E \to I_Z(3) \to 0 \ \ ,
$$
with $Z$ a smooth codimension $2$ subvariety of $\P^n$. In our hypothesis, $Z$ is
contained in a hyperplane, and $I_Z$ is generated by cubic forms so that $Z$ is
either a linear subspace, a complete intersection of type $(1,2)$ or of type $(1,3)$.
When $Z$ is a linear subspace $\Lambda _Z \cong \Oc(-1)$ is not globally generated.
When $Z$ is a complete intersection of type $(1,2)$ $\Lambda _Z \cong \Oc _Z$. Then 
$E \cong \Oc (1) \oplus \Oc (2)$. When $Z$ is a complete intersection of type 
$(1,3)$, $\Lambda _Z \cong \Oc (1)$ is generated by $(n-1)$ sections, so that $E$ is
given by an extension:
$$0 \to (n-1)\Oc \to E \to I_Z(3) \to 0 \ \ .
$$
Then $E$ is given by an extension:
$$
0 \to \Oc(-1) \to n\Oc \oplus \Oc (2) \to E \to 0 \ \ .
$$

II. {\bf Assume} $H^0(E(-2))=0$ and $H^0(E(-1))\not =0$. 

Then in the extension (***) $Z$ is not contained in a hyperplane, is contained in at 
least one quadric and $I_Z$ is generated by cubics. Then $3 \le \deg (Z)\le 6$.

(a) $\deg (Z) =3$. According to \cite{W}, $Z$ is one of: 

(a1) $Z \subset \P^5$, image of a Segre embedding of $\P^1 \times P^2$ in $P^5$, 

(a2) $Z\subset \P^4$, a smooth hyperplane section of the first case. 

The case (a1) does not give an indecomposable vecor bundle. Indeed, $I_Z$ would 
have a resolution of the shape:
$$
0 \to 2 \Oc_{\P^5} (-3) \to 3 \Oc _{\P^5}(-2) \to I_Z \to 0 \ \ ,
$$
(cf. \cite{BaMa} ) and $\Lambda _Z |_{\P^1 \times \P^2} \cong \Oc_{\P^1 
\times \P^2}(1,0) $. It is easy to see that, in this case, $E \cong 3\Oc (1)$.

In the case (a2) one gets for $Z'=Z\cap \P^4$ a similar resolution:
$$
0 \to 2 \Oc_{\P^4} (-3) \to 3 \Oc _{\P^4}(-2) \to I_{Z'} \to 0 \ \ ,
$$
and again $E$ is direct sum of line bundles.

(b) $\deg(Z)=4$. The only codimension $2$, smooth projective varieties of degree 
$4$, not contained in a hyperplane, contained in a hyperquadric and with the ideal 
generated by $H^0(I_Z(3))$ are the complete intersections of type $(2,2)$, (cf.
\cite{S-D}, or \cite{Ok}). In this case one gets a resolution for $E$:
$$
0 \to \Oc(-1) \to (n-1)\Oc \oplus 2 \Oc(1) \to E \to 0 \ \ .
$$

(c) $\deg (Z) =5$. The only codimension $2$, smooth projective varieties of degree $5$,
not contained in a hyperplane, with $H^0(I_Z(2)) \not =0$ and with the ideal generated
by degree $3$ component are (cf. \cite{Ok}, \cite{Ok-bis} or \cite{Io}):

(c1) Castelnuovo surface $Z \cong \P^2(x_0, \ldots ,x_7)$, i.e. $P^2$ blown-up in 8 points in general 
position, embedded in $\P^4$ via the quartics through $x_i$, with multiplicity $2$ 
in $x_0$. $I_Z$ has a minimal resolution of the shape (see \cite{Ok}, proof of 2.9):
$$
0 \to 2\Oc (-4) \to 2\Oc(-3)\oplus \Oc (-2) \to I_Z \to 0 \ \ .
$$
$\Lambda _Z \cong \omega _Z(2)$ is generated by global sections. As 
$h^0(\omega _Z(3))=8$, one gets a vector bundle $F$ of rank $9$ with a resolution 
of the shape:
$$
0 \to 2\Oc (-1) \to 10\Oc \oplus \Oc (1) \to F \to 0 \ \ .
$$
$F$ is an extension by a trivial bundle of a vector bundle $E$ of $\rank =4$. 
Indeed, by the property {\bf A.}, $F$ contains a trivial subbundle of rank $5$. 
Let $E$ be the quotient. By standard arguments one sees 
that this lifts to a trivial subbundle of $10\Oc\oplus \Oc(1)$ such that one has 
the commutative diagram;

$$
\begin{CD}
@. @. 0 @. 0 @. \\
@. @. @VVV @VVV @.\\
@. @. 5\Oc @= 5\Oc @. \\
@. @. @VVV @VVV @.\\
0 @>>> 2\Oc(-1) @>>>10\Oc \oplus \Oc(1) @>>> F @>>> 0 \\
@. @| @VVV @VVV @.\\
0 @>>> 2\Oc(-1) @>>>5\Oc \oplus \Oc(1)  @>>> E @>>> 0 \\   
@. @. @VVV @VVV @.\\   
  @. @. 0 @. 0 @.      
\end{CD}
$$

(c2) $Z$ = Castelnuovo threefold in $\P^5$. In this case the resolution of $I_Z$ has 
the shape (cf. \cite{Ok-bis}, proof of 1.5):
$$
0 \to 2\Oc (-4) \to 2\Oc(-3)\oplus \Oc (-2) \to I_Z \to 0 \ \ .
$$

Like above, one obtains a vector bundle $E$ on $\P^5$ with a resolution of the form:
$$
0 \to 2\Oc (-1) \to 6\Oc \oplus \Oc (1) \to E \to 0 \ \ .
$$

(d) $\deg (Z)=6$. The condition that $Z$ is contained in a quadric makes possible only:
$Z$ = complete intersection of type $(2,3)$. By standard arguments one obtains
a vector bunlde $E$ on $\P^n$ of $\rank =n$ with a resolution  of the shape:
$$
0 \to \Oc(-2) \to n\Oc \oplus \Oc(1) \to E \to 0 \ \ .
$$

III. {\bf Assume } $H^0(E(-1))=0$. We have to consider the smooth projective 
varieties of codimension $2$, not contained in any hyperquadric, but with ideal
generated by cubic forms.

(a) $\deg (Z) =4$. The only possibility is the Veronese surface in $\P^4$.
It is known (cf. \cite{Ok}, Corr. 2.8) that the ideal of $Z$ has an 
$\Omega$-resolution of the form:
$$
0 \to 3\Oc \to \Omega ^1(2) \to I_Z(3) \to 0 \ \ .
$$
In fact any vector bundle $E$ in an extension:
$$
0 \to 3\Oc \to E \to I_Z(3) \to 0 \ \ .
$$
with $Z$ the Veronese surface in $\P^4$ is isomorphic to $\Omega^1(2)$.
This can be seen via the Beilinson spectral sequence (cf. \cite{OSS}, II, 3.1.3)
associated to $E(-1)$, which turns out to be (very) degenerated.

(b) $\deg (Z)=5$. In our conditions $Z$ is an elliptic scroll with $e=-1$ and
$H \equiv C +2f$ (cf. \cite{Ok}, Theorem 2.10). In \cite{Ok} it is shown that 
$I_Z$ has a resolution:
$$
0 \to T_{\P^4}(-2) \to 5\Oc \to I_Z(3) \to 0
$$
From here one gets a $\rank$ $3$ vector bundle as an extension:
$$
0 \to 2\Oc \to E \to I_Z(3) \to 0
$$
The resolution of $I_Z(3)$ provides a resolution of $E$:

$$
0 \to T_{\P^4}(-2) \to 7\Oc \to E \to 0 \ \ ,
$$
which shows that $E$ is the dual of a Trautmann-Vetter (Tango) bundle on $\P^4$.
(cf. \cite{T}, \cite{V}, \cite{Tango}, \cite{KPR} ).

(c) $\deg(Z)=6$. The only possibilities are the Bordiga surface in $\P^4$ and the
Bordiga threefold in $\P^5$ (cf. \cite{Ok-bis} or \cite{Io}).

In both cases the minimal resolution has the shape:
$$
0 \to 3\Oc (-4) \to 4\Oc(-3) \to I_Z \to 0 \ \ .
$$
These two cases will produce a $\rank$ $4$ bundle on $\P^4$ with a resolution:
$$
0 \to 3\Oc_{\P^4}(-1) \to 7\Oc_{\P^4}\to E \to 0
$$
and a $\rank$ $5$ bundle on $\P^5$ with a resolution:
$$
0 \to 3\Oc_{\P^5}(-1) \to 8\Oc_{\P^5}\to E \to 0\ \ .
$$

(d) $\deg(Z)=7$. We could use \cite{Ok-tert}, but one can proceed directly.
Indeed, in our situation $Z$ is linked in a complete intersection of type $(3,3)$ to
a locally Cohen-Macaulay subscheme $Y$ of $\P^n$, $(n=4,5)$, of pure codimension $2$ 
and degree $2$.  $Y$ cannot be a double structure on a linear subspace of codimension
$2$, as such a structure does not exist (cf. \cite{Ma2}) and cannot be a union of
two linear subspaces of codimension $2$ meeting along a linear subsapce of 
codimension $4$, as this would contradict Cohen-Macaulay property. Then $Y$
is a quadric of codimension $2$, degenerated or not. It follows the resolution of
$I_Z$:
$$
0 \to \Oc(-5)\oplus \Oc(-4) \to 3\Oc(-3) \to I_Z \to 0
$$
These $Z$'s give rise to the vector bundles on $\P^4$, respectively $\P^5$:
$$
0 \to \Oc_{\P^4}(-2)\oplus \Oc_{\P^4}(-1) \to 6\Oc_{\P^4} \to E \to 0 \ \ ,
$$
$$
0 \to \Oc_{\P^5}(-2)\oplus \Oc_{\P^5}(-1) \to 7\Oc_{\P^5} \to E \to 0 \ \ ,
$$
of ranks $4$, respectively $5$.

(e) $\deg(Z)=8$. This situation does not occur in our context, as such a $Z$, being 
linked to linear subspace would have $I_Z(3)$ not globally generated.

(f) $\deg(Z)=9$. Then $Z$ is a complete intersection and the vector bundles obtained 
this way are:
$$
0 \to \Oc_{\P^4}(-3) \to 5\Oc_{\P^4} \to E \to 0 \ \ ,
$$
of rank $4$ on $\P^4$ and
$$
0 \to \Oc_{\P^5}(-3) \to 6\Oc_{\P^5} \to E \to 0 \ \ ,
$$
of rank $5$ on $\P^5$.

\section{Vector Bundles on $\P^n$, $n\ge 5$}

From  (iii) and (iv) in the Theorem it follows directly (v) for $\P^5$.
We prove (v) by induction by standard (known) cohomological arguments.

\begin{lem}
If $E$ is a vector bundle on $\P^{n+1}$ and its restriction $E'$ to a hyperplane 
has $H^i(E'(j))=0$ for $1 \le i \le n-2$ and any $j$, then $H^i(E(j))=0$ for
$1 \le i \le n-1$ and any $j$.
\end{lem}

\proof From the exact sequence:
$$
0 \to E(-1) \to E \to E' \to 0
$$
it follows, for all $j$:
\begin{eqnarray*}
H^i(E(j-1)) \cong H^i(E(j)), & \ \ 2\le i \le n-2 \ \ , \\
H^1(E(j-1)) \to H^1(E(j))  & \hbox{ surjective} \\
H^{n-1}(E(j-1)) \to H^{n-1}(E(j)) & \hbox{ injective}
\end{eqnarray*}
As $H^i(E(j))=0$ for $i =1, \ldots n$ and $j \ll 0$ or $j \gg 0$ the assertion follows.
\qed
\begin{lem}
If $E$ is a vector bundle on $\P^n$ with the property $H^i(E(j))=0$ for $1 \le i \le n-2$
then it has a resolution of length $1$:
$$
0 \to L_1 \to L_0 \to E \to 0 \ \ ,
$$
where $L_0,\ L_1$ are direct sums of line bundles.
\end{lem}

\proof Indeed, take a (minimal) surjection $L_0 \to E \to 0$, with $L_0$  a direct 
sum of line bundles and denote by $L_1$ the 
kernel. It follows $H^i((L_1(j))=0$ for $i=1, \ldots n-1$ and all $j$. By the spliting 
criterion of Horrocks (cf. \cite{Ho}, or \cite{OSS} I, Theorem  2.3.1) it follows that 
$L_1$ is a direct sum of line bundles.
\qed

The assertion (v) in the Theorem follows directly from the above lemmas.

\end{document}